\documentclass{amsart}
\usepackage{amsfonts}
\usepackage{amssymb,latexsym}
\usepackage{stmaryrd}
\usepackage{amsmath}  %Overlap with newcommand \mod
\usepackage{verbatim}

\input xy  
\xyoption{all}
\xyoption{frame}

\theoremstyle{plain} % style plain
\newtheorem{theorem}{Theorem}[section]
\newtheorem{lemma}[theorem]{Lemma}

\newtheorem{proposition}[theorem]{Proposition}
\newtheorem{corollary}[theorem]{Corollary}

\theoremstyle{definition} % style definition
\newtheorem{example}[theorem]{Example}

\newtheorem{remark}[theorem]{Remark}

\newcommand{\Ker}{\mbox{\rm Ker\,}}

\newcommand{\End}[1]{\operatorname{\rm End}_{#1}}

\newcommand{\Hom}[1]{\operatorname{{\rm Hom}}_{#1}}

\newcommand{\rad}{\operatorname{{\rm rad}}}
\newcommand{\soc}{\operatorname{{\rm soc}}}

\newcommand{\Spec}{\operatorname{\rm Spec\,}}

\newcommand{\pres}[1]{\mbox{{\rm pres}}_{#1}}

\newcommand{\add}{\mbox{{\rm add \!}}}

\newcommand{\MOD}{\mbox{{\rm mod \!}}}

\newcommand{\rep}{\mbox{{\rm rep \!\!}}}

\newcommand{\proj}{\operatorname{{\rm proj }}}

\newcommand{\demo}[1]{\textsc{Proof.} #1 \hfill $\Box$ \bigskip}

\newcommand{\cA}{\mathcal{A}}

\newcommand{\cC}{\mathcal{C}}
\newcommand{\cD}{\mathcal{D}}

\newcommand{\cH}{\mathcal{H}}

\newcommand{\cO}{\mathcal{O}}

\newcommand{\cY}{\mathcal{Y}}
\newcommand{\cZ}{\mathcal{Z}}
\newcommand{\bC}{\mathbb{C}}
\newcommand{\bN}{\mathbb{N}}

\newcommand{\bZ}{\mathbb{Z}}
\newcommand{\bA}{\mathbb{A}}
\newcommand{\bD}{\mathbb{D}}
\newcommand{\bE}{\mathbb{E}}
\newcommand{\bd}{\mathbf{d}}

\newcommand{\fM}{\mathfrak{M}}
\newcommand{\spl}{\mathrm{split}}

\newcommand{\Si}{\Sigma}
\newcommand{\g}{\mathfrak{g}}
\newcommand{\eg}{\emph{e.g. }}

\begin{document}

\title[Nakajima varieties and repetitive algebras]{Nakajima varieties and repetitive algebras}
\author[B. Leclerc \and P.-G. Plamondon]{Bernard Leclerc \and Pierre-Guy Plamondon}

\address{ Bernard Leclerc \\
          Normandie Univ. France\\
          UNICAEN, LMNO F-14032 Caen, France\\
          CNRS UMR 6139, F-14032 Caen, France\\
          Institut Universitaire de France.}
\address{ Pierre-Guy Plamondin \\
Normandie Univ. France\\
          UNICAEN, LMNO F-14032 Caen, France\\
          CNRS UMR 6139, F-14032 Caen, France}

          \email{bernard.leclerc@unicaen.fr}
\email{pierre-guy.plamondon@unicaen.fr}

\thanks{Communicated by H. Nakajima. Received September 22, 2012.
Revised January 18, 2013.
The second author was financially supported by the FQRNT}
%\date{\today}

%\classification{18E30}

\begin{abstract}
We realize certain graded Nakajima varieties of finite Dynkin type 
as orbit closures of repetitive algebras of Dynkin quivers.
As an application, we obtain that the perverse sheaves introduced
by Nakajima for describing irreducible characters of quantum loop 
algebras are isomorphic to the intersection cohomology sheaves
of these orbit closures.  
\end{abstract}

\maketitle

\tableofcontents

%--------------------------------------------------------------------------
\section{Introduction}
%..........................................................................

Let $U_q(L\g)$ be the quantum loop algebra of a simple Lie algebra $\g$ 
of Dynkin type $\bA, \bD$ or $\bE$, introduced by Drinfeld and Jimbo. This is 
a $\bC$-algebra, and we assume that the quantum parameter $q\in\bC^*$
is not a root of unity. The category of finite-dimensional representations
of $U_q(L\g)$ has been studied by many authors \cite{AK, CP1, FR, GV, H1, NJAMS}. 

In type $\bA$, a geometric construction of the standard and irreducible 
$U_q(L\g)$-modules was given by Ginzburg and Vasserot \cite{GV,V}. 
In particular, they showed that the composition multiplicities of the 
standard modules can be expressed in terms of the local intersection cohomology of graded nilpotent 
orbit closures. This theory was extended by Nakajima to all 
simply-laced types. To do so, he replaced the graded nilpotent orbit closures
by certain graded versions $\fM_0^\bullet(V,W)$ 
of the quiver varieties $\fM_0(V,W)$ previously introduced in his
geometric construction of the irreducible finite-dimensional $\g$-modules.
Nakajima's varieties are defined as quotients in geometric invariant theory,
and in general they do not have any other more explicit description.
 
%These varieties depend on some graded vector spaces $V$ and $W$ encoding
%respectively a weight 

Let $Q$ be a Dynkin quiver of the same type as $\g$.
It was recently observed \cite{HL11} that the orbit closures of $Q$
are isomorphic to some particular Nakajima's graded quiver varieties
$\fM_0^\bullet(V,W)$.
In this paper we generalize this observation and show that 
all the orbit closures of the repetitive algebra \cite{HW83} of the path algebra $\bC Q$ 
are isomorphic to certain varieties $\fM_0^\bullet(V,W)$.
The repetitive algebra $\widehat{A}$ of an algebra~$A$, 
first introduced by Hughes and Waschb\"usch,
is a selfinjective infinite-dimensional algebra,
defined as a kind of infinite matrix algebra 
(its precise definition is recalled in \S\ref{sect3.1}).
It was shown by Happel that, if $A$ has finite global dimension, 
the stable module category of $\widehat{A}$
is equivalent to the derived category of $A$.  
When $Q$ is a linearly oriented Dynkin quiver of type $\bA$, the
orbit closures of the repetitive algebra~$\widehat{\bC Q}$
coincide with the graded nilpotent orbit closures of Ginzburg and Vasserot.
Our result may thus be regarded as an explicit description, similar
to that of Ginzburg and Vasserot, of a large collection of graded
Nakajima varieties $\fM_0^\bullet(V,W)$.

Although we only realize in this way the varieties corresponding to 
some particular $W$'s, we get enough of them to express all the
irreducible characters of $U_q(L\g)$ in terms of the local intersection
cohomology of the orbit closures of $\widehat{\bC Q}$.
Recall that the perverse sheaves used by Lusztig \cite{L} in his geometric construction
of the canonical basis of $U_q^+(\g)$ are precisely the intersection cohomology
sheaves of the orbit closures of $Q$.
This description of the canonical basis was inspired by Ringel's theorem \cite{R}, stating that
$U_q^+(\g)$ is isomorphic to the Hall algebra of $Q$.
Similarly our result shows that the perverse sheaves used by Nakajima
in his geometric construction of the canonical basis of the Grothendieck 
ring of $U_q(L\g)$ are essentially the intersection cohomology
sheaves of the orbit closures of $\widehat{\bC Q}$. 
This was inspired by the recent result of \cite{HL11}, stating that
the $t$-deformed Grothendieck ring is isomorphic to T\"oen's Hall
algebra of the derived category of~$Q$, which, by Happel's theorem,
is isomorphic to T\"oen's algebra of the stable category of $\widehat{\bC Q}$.

Another interesting application of our results is that we can apply
Nakajima's well-developed theory to get new information on the orbit closures of the repetitive
algebras of Dynkin type. For instance we obtain desingularizations
of these orbit closures with favorable properties, like the vanishing of odd cohomology
groups of the fibers. Also, using Nakajima's algorithm for computing the
irreducible characters of $U_q(L\g)$ \cite{N}, we can in principle calculate the local intersection
cohomology of these orbit closures.

Our paper is structured as follows. 
In \S\ref{sect2}, we recall the definition of the graded Nakajima varieties
$\fM_0^{\bullet}(V,W)$ and $\fM_0^{\bullet}(W)$, as well as the stratification of 
$\fM_0^{\bullet}(W)$ by the strata $\fM_0^{\bullet \mathrm{reg}}(V,W)$.
We then show that the variety $\fM_0^{\bullet}(W)$ is isomorphic to the 
variety of representations with dimension vector $\bd$ of a certain finite-dimensional
algebra $e_\Si\Lambda e_\Si$ obtained by projectivization (Theorem~\ref{theo::iso of varieties}).  
Here $\bd$ is the graded dimension of $W$. 
If $W$ is arbitrary, $e_\Si\Lambda e_\Si$ looks rather complicated.
But we show in \S\ref{sect3} that for certain special choices of $W$,
the algebra $e_\Si\Lambda e_\Si$ is isomorphic to the repetitive algebra
$\widehat{A}$ of $A=\bC Q$ (Lemma~\ref{lemm::projectives}). Moreover, in this case,
we show (Theorem~\ref{main_theorem}) that Nakajima's stratification of $\fM_0^{\bullet}(W)$ 
is identical to the usual orbit stratification of $\rep_{\bd}(\widehat{A})$ (recall that 
since $Q$ is of Dynkin type, $\widehat{A}$ is locally of finite representation
type). The precise correspondence between isoclasses $N$ of $\widehat{A}$-modules
of dimension $\bd$ and strata $\fM_0^{\bullet \mathrm{reg}}(V,W)$ is described
in Corollary~\ref{coro::explicit bijection}.    
Finally in \S\ref{sect4}, we apply our results to the representation theory of
quantum loop algebras $U_q(L\g)$. 
We start by recalling the definition of the tensor category
$\cC_\bZ$ of finite dimensional $U_q(L\g)$-modules, and we
state Nakajima's theorem giving a geometric description of 
the composition multiplicities
of the standard modules of $\cC_\bZ$.
Then we deduce from Theorem~\ref{main_theorem} that there is a 
natural one-to-one correspondence between the isomorphism classes
of objects in the stable category $\underline{\MOD} \widehat{A}$
and the isomorphism classes of simple objects in $\cC_\bZ$
(Proposition~\ref{parametrization}). We conclude with the description
of the composition multiplicities of the standard modules in terms
of the intersection cohomology of the $\widehat{A}$-orbit closures.

%--------------------------------------------------------------------------
\section{Nakajima varieties and representations of algebras}\label{sect2}
%...........................................................
\subsection{Varieties of representations}\label{sect2.1}
We will work over the field $\bC$ of complex numbers.  
Let $Q=(Q_0,Q_1)$ be a quiver with vertex set $Q_0$ and arrow set $Q_1$.  Let $s,t:Q_1 \rightarrow Q_0$ be maps taking each arrow $a$ to its source $s(a)$ or its target $t(a)$.
The number of vertices may be infinite, but we will assume that the number of arrows 
incident with any given vertex is finite.  
Let $A := \bC Q/ R$ be the quotient of the path algebra of $Q$ by an ideal $R$.  We will always assume that $R$ is contained in the two-sided ideal generated by the arrows of $Q$.  Then $A$ is an associative algebra which has no unit if the number of vertices of $Q$ is infinite.  For any vertex $i$ of $Q$, let $e_i$ be the trivial path at $i$.  

A \emph{representation} of $(Q,R)$ is a pair $\big( (V_i)_{i\in Q_0}, (\varphi_a)_{a\in Q_1} \big)$, where each $V_i$ is a vector space over $\bC$ and each $\varphi_a$ is a linear map from $V_{s(a)}$ to $V_{t(a)}$ such that all polynomials in the $\varphi_a$ corresponding to elements in $R$ vanish.  The representations of $(Q,R)$ form a category.  This category is equivalent to the category of left modules $M$ over $\bC Q$ which have the property that $M = \bigoplus_{i\in Q_0} e_i M$.

Let $\bd \in \bN^{Q_0}$ be a finite dimension vector, that is, such that the sum of its entries $d_i$ for $i\in Q_0$ is finite.  The \emph{variety of representations of $A$ of dimension vector $\bd$} is the closed subset $\rep_\bd (A)$  of
\begin{displaymath}
	\bigoplus_{a\in Q_1} \Hom{\bC}(\bC^{d_{s(a)}}, \bC^{d_{t(a)}})
\end{displaymath}
which vanishes on polynomials defined by the relations $R$.  Its points are in bijection with the representations of $A$ with dimension vector $\bd$ and underlying vector space $\prod_{i\in Q_0}\bC^{d_i}$.  It is a finite-dimensional affine algebraic variety.

The affine algebraic group
\begin{displaymath}
	G_\bd := \prod_{i\in Q_0} GL(\bC^{d_i})
\end{displaymath}
acts on $\rep_\bd(A)$ by base change at every vertex.  The $G_{\bd}$-orbits are in bijection with the isomorphism classes of $A$-modules with dimension vector $\bd$.

%............................................................
\subsection{Nakajima's varieties} \label{sect2.2}
Let $Q$ be a quiver of Dynkin type $\bA$, $\bD$ or~$\bE$.  
For $i, j$ in $Q_0$, we write $i\sim j$ if there is in $Q_1$ an arrow 
$i\to j$ or an arrow $j\to i$.
A \emph{height function} on $Q$ is a map
\begin{displaymath}
	\xi : Q_0 \longrightarrow \bZ
\end{displaymath}
such that for any arrow $i\rightarrow j$ in $Q$, we have that $\xi_j = \xi_i -1$.  Since $Q$ is a tree, there exists a height function on it.

Define the \emph{repetition quiver} $\widehat{Q}$ thus:
\begin{itemize}
	\item the set of vertices of $\widehat{Q}$ is $\widehat{Q}_0 = \{ (i,p) \in Q_0\times \bZ \ | \ p - \xi_i \in 2\bZ \}$;  
	\item for any arrow $c:i\rightarrow j$ in $Q$ and any $(i,n) \in \widehat{Q}_0$, there is an arrow $(c,n):(i,n)\rightarrow (j,n-1)$ and an arrow $(\overline{c},n-1):(j,n-1)\rightarrow (i,n-2)$ in~$\widehat{Q}$.
\end{itemize}

Now, let $\widehat{\Gamma}$ be the quiver which contains $\widehat{Q}$ as a full subquiver, and which has additionnal vertices and arrows as follows:
\begin{itemize}
	\item one additionnal vertex for each element of $\{ (i,p) \in Q_0\times \bZ \ | \ p - \xi_i \in (1 + 2\bZ) \}$. In particular, the set of vertices of $\widehat{\Gamma}$ is $Q_0 \times \bZ$;
	\item additionnal arrows $a_{(i,n+1)}:(i, n+1) \rightarrow (i,n)$ and $b_{(i,n)}:(i,n) \rightarrow (i, n-1)$ for every $(i,n) \in \widehat{Q}_0$. 
\end{itemize}
We will study representations of $\widehat{\Gamma}$ satisfying the following set of relations:
\begin{displaymath}
	a_{(i,n-1)}b_{(i,n)} + \sum_{c:i\rightarrow j}(\overline{c},n-1)(c,n) - \sum_{c:j\rightarrow i}(c,n-1)(\overline{c},n) = 0, \ (i,n)\in \widehat{Q}_0.
\end{displaymath}
Let $R$ be the ideal generated by these relations.  Let $\Lambda$ be the quotient of $\bC\widehat{\Gamma}$ by $R$.

Let $\bd$ be a finite dimension vector on $\widehat{\Gamma}$.  We will use the following notation:
\begin{equation}\label{eqVW}
	V = \bigoplus_{(i,n)\in \widehat{Q}_0} \bC^{d_{(i,n)}}, 
\quad W = \bigoplus_{(i,n)\in \widehat{\Gamma}_0 \setminus \widehat{Q}_0} \bC^{d_{(i, n)}},
\end{equation}
where $V$ is a $\widehat{Q}_0$-graded vector space, and $W$ is a $(\widehat{\Gamma}_0 \setminus \widehat{Q}_0)$-graded vector space.
Then the variety of representations $\rep_{\bd}(\Lambda)$ is a closed subset of the direct sum of these three vector spaces:
\begin{displaymath}
	L^{\bullet}(V,W) = \bigoplus_{(i,n)\in \widehat{Q}_0} \Hom{\bC}(V_{(i,n)}, W_{(i,n-1)});
\end{displaymath}
\begin{displaymath}
	L^{\bullet}(W,V) = \bigoplus_{(i,n)\in \widehat{Q}_0} \Hom{\bC}(W_{(i,n+1)}, V_{(i,n)}); \\
\end{displaymath}
\begin{displaymath}
	E^{\bullet}(V) = \bigoplus_{i\rightarrow j, (i,n)\in \widehat{Q}_0} \Hom{\bC}(V_{(i,n)}, V_{(j,n-1)}) \oplus \bigoplus_{i\rightarrow j, (j,m)\in \widehat{Q}_0} \Hom{\bC}(V_{(j,m)}, V_{(i,m-1)}).
\end{displaymath}
The variety $\rep_{\bd}(\Lambda)$ will also be denoted by 
$\Lambda^{\bullet}(V,W)$ when we want to emphasize this decomposition.

The affine algebraic group
\begin{displaymath}
	G_V := \prod_{(i,n)\in \widehat{Q}_0}GL(V_{(i,n)})
\end{displaymath}
acts naturally on $\Lambda^{\bullet}(V,W)$ by base change at vertices $(i,n)\in \widehat{Q}_0$.  Define the affine quotient
\begin{displaymath}
	\fM_0^{\bullet}(V,W) := \Lambda^{\bullet}(V,W) \sslash G_V := \Spec \bC [\Lambda^{\bullet}(V,W)]^{G_V}.
\end{displaymath}
If $V \subset V'$ as $\widehat{Q}_0$-graded vector spaces, then there is a natural closed embedding $\fM_0^{\bullet}(V,W) \subset \fM_0^{\bullet}(V',W)$.  Define
\begin{displaymath}
	\fM_0^{\bullet}(W) = \bigcup_V \fM_0^{\bullet}(V,W).
\end{displaymath}
It is an affine variety on which the algebraic group
\begin{displaymath}
	G_W := \prod_{(i,n)\in \widehat{\Gamma}_0 \setminus\widehat{Q}_0} GL(W_{(i, n)})
\end{displaymath}
acts by base change at vertices $(i,n)\in \widehat{\Gamma}_0 \setminus \widehat{Q}_0$.  
Nakajima proves (see for instance \cite[Section 3]{Nak11}) that there is a stratification
\begin{displaymath}
	\fM_0^{\bullet}(W) = \bigsqcup_V \fM_0^{\bullet \mathrm{reg}}(V,W),
\end{displaymath}
where $\fM_0^{\bullet \mathrm{reg}}(V,W)$ is the open subset of $\fM_0^{\bullet}(V,W)$ parametrizing the closed free $G_V$-orbits.  He also proves that a necessary condition for $\fM_0^{\bullet \mathrm{reg}}(V,W)$ to be non-empty is that
\begin{displaymath}
	\dim W_{(i,n)} - \dim V_{(i,n+1)} - \dim V_{(i,n-1)} + \sum_{j \sim i} \dim V_{(j,n)} \geq 0,
\end{displaymath}
for all $(i,n) \in \Gamma_0 \setminus Q_0$.  In that case, the pair $(V,W)$ is called \emph{dominant}.

Similarly, we will say that a representation $N$ of $\Lambda$ is \emph{dominant} if the associated pair $(V,W)$ is dominant.  The representation $N$ is \emph{stable} if any graded subspace of $V$ stable under the action of all $(c,n)$ and $(\overline{c},n+1)$ and vanishing under the action of all $b_{(i,n)}$ is the zero subspace.  In other words, $N$ is stable if it has no non-zero subrepresentations supported on $V$.

The stratum $\fM_0^{\bullet \mathrm{reg}}(V,W)$ is non-empty if, and only if, there exists a stable representation $N$ whose underlying graded vector space is $(V,W)$, and $(V,W)$ is dominant.

\begin{example}
Consider the quiver $Q$ of type $\bA_4$ given by
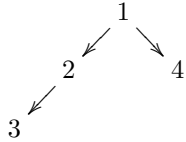
\begin{displaymath}
	\xymatrix@-1.2pc{ & & 1\ar[dr]\ar[dl] & \\
	           & 2\ar[dl] & & 4 \\
	           3 & & &
	}
\end{displaymath}
with height function defined by $\xi_1 = 2$, $\xi_2 = \xi_4 = 1$ and $\xi_3 = 0$.  Then the points of $\Lambda^{\bullet}(V,W)$ are given by collections of morphisms which can be organized as in the following picture:
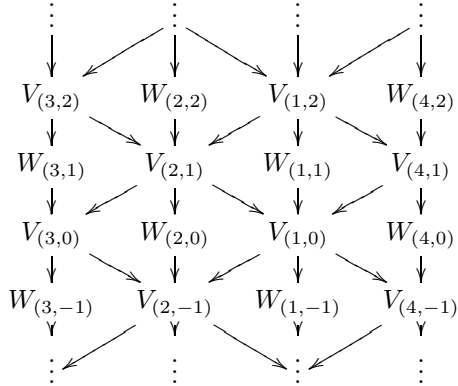
\begin{displaymath}
	\xymatrix@-1.2pc{ \vdots\ar[d] & \vdots\ar[dl]\ar[d]\ar[dr]  & \vdots\ar[d]  & \vdots\ar[dl]\ar[d] \\
	                  V_{(3,2)}\ar[dr]\ar[d] & W_{(2,2)}\ar[d] & V_{(1,2)}\ar[dl]\ar[d]\ar[dr] & W_{(4,2)}\ar[d] \\
	                  W_{(3,1)}\ar[d] & V_{(2,1)}\ar[dr]\ar[d]\ar[dl] & W_{(1,1)}\ar[d] & V_{(4,1)}\ar[dl]\ar[d] \\
	                  V_{(3,0)}\ar[dr]\ar[d] & W_{(2,0)}\ar[d] & V_{(1,0)}\ar[dl]\ar[d]\ar[dr] & W_{(4,0)}\ar[d] \\
	                  W_{(3,-1)}\ar[d] & V_{(2,-1)}\ar[dr]\ar[d]\ar[dl] & W_{(1,-1)}\ar[d] & V_{(4,-1)}\ar[dl]\ar[d] \\
	                 \vdots & \vdots  & \vdots  & \vdots  
	}
\end{displaymath}
The quiver $Q$ is found in this picture by taking the full subquiver whose vertices are $V_{(1,2)}$, $V_{(2,1)}$, $V_{(3,0)}$ and $V_{(4,1)}$.  The quiver whose arrows and vertices are all the $V$'s and $W$'s is $\widehat{\Gamma}$; the full subquiver whose vertices are the $V$'s is $\widehat{Q}$. 
\end{example}

%.............................................................
\subsection{Projectivization}

Let $\widehat{\Gamma}$, $R$ and $\Lambda$ be as above.  
Let $\Si$ be a subset of~$\widehat{\Gamma}_0$.  
Define $\widehat{\Gamma}_\Si$ as the full subquiver of $\widehat{\Gamma}$  
whose vertices are those $x$ in $\widehat{\Gamma}_0$
such that there exist two vertices $w$ and $y$ in $\Si$ and a path from $w$ to $y$ 
passing through~$x$. (One may think of $\widehat{\Gamma}_\Si$ as a kind of convex hull of $\Si$).

Let $R_\Si := R \cap \bC \widehat{\Gamma}_\Si$, and define
$  \Lambda_\Si := \bC \widehat{\Gamma}_\Si /R_\Si$.
If $\bd$ is a dimension vector supported on $\widehat{\Gamma}_\Si$, then
$\rep_{\bd}(\Lambda) \cong \rep_{\bd}(\Lambda_{\Si})$.
Note that if $\Si$ is finite, then $\widehat{\Gamma}_\Si$ only has a finite number of vertices,
and $\Lambda_\Si$ is a finite-dimensional algebra with a unit.

For the rest of the section, we will assume that $\Si$ is finite, and that 
\begin{equation}\label{eqn::assumption}
\Si \subset \widehat{\Gamma}_0 \setminus \widehat{Q}_0.  
\end{equation}
Define the idempotent
\[
e_\Si := \sum_{x \in \Si} e_x,
\]
and consider the algebra 
$
e_\Si\Lambda e_\Si = e_\Si\Lambda_\Si e_\Si.
$
We will list here some facts concerning this algebra, taken mainly from Chapter II.2 of \cite{ARS}.

The algebra $e_\Si\Lambda_\Si e_\Si$ is a finite-dimensional algebra with unit $e_\Si$, and the $e_x$, for $x$ in 
$\Si$, form a complete set of pairwise orthogonal primitive idempotents.  
It is isomorphic as an algebra to $\End{\Lambda_\Si}(\Lambda_\Si e_\Si)^{\mathrm{op}}$.  
From this isomorphism, we see that we have a functor
\begin{displaymath}
	\Hom{\Lambda_\Si}(\Lambda_\Si e_\Si, ?) : \MOD \Lambda_\Si \longrightarrow \MOD e_\Si \Lambda_\Si e_\Si,
\end{displaymath}
which is exact since $\Lambda_\Si e_\Si$ is a projective $\Lambda_{\Si}$-module.  The process of applying this functor is called \emph{projectivization} in \cite{ARS}.

Let $\pres{} (\Lambda_\Si e_\Si)$ be the full subcategory of $\MOD \Lambda_\Si$ 
consisting of modules $M$ which admit a projective presentation
\begin{displaymath}
	P_1 \longrightarrow P_0 \longrightarrow M \longrightarrow 0,
\end{displaymath}
with $P_1$ and $P_0$ in $\add (\Lambda_\Si e_\Si)$.  

\begin{proposition}[Proposition II.2.5 of \cite{ARS}]
The functor $\Hom{\Lambda_\Si}(\Lambda_\Si e_\Si, ?)$ restricts to an equivalence of categories
\begin{displaymath}
	\pres{} (\Lambda_\Si e_\Si) \longrightarrow \MOD e_\Si \Lambda_\Si e_\Si.
\end{displaymath}
\end{proposition}

Let $M$ be a $\Lambda_\Si$-module.  
Then the dimension vector of its image in $\MOD e_\Si \Lambda_\Si e_\Si$ is obtained from $\dim M$ by forgetting the components which are not in $\Si$.  This is seen from the isomorphisms
\begin{displaymath}
	e_x \Hom{\Lambda_\Si}(\Lambda_\Si e_\Si, M) \cong e_x \Hom{\Lambda_\Si}(\Lambda_\Si , M) \cong e_x M.
\end{displaymath}

For any dimension vector $\bd$ on $\widehat{\Gamma}$, we write $\bd = \bd_V + \bd_W$, where $\bd_V$ is supported on $\widehat{Q}_0$ and $\bd_W$ is supported on $\widehat{\Gamma}_0 \setminus \widehat{Q}_0$.  Then the above discussion gives us the following, which we will use later.

\begin{corollary}\label{coro::functor}
Let $M$ be a module over $e_\Si \Lambda_\Si e_\Si$.  Then there exists a $\Lambda_\Si$-module $L$ such that $\Hom{\Lambda_\Si}(\Lambda_\Si e_\Si, L)$ is isomorphic to $M$.  Its dimension vector $\bd = \bd_V + \bd_W$ is such that $\bd_W$ is the dimension vector of $M$.
\end{corollary}

%--------------------------------------------------------------------------
\subsection{An isomorphism of varieties}

In this section, we prove the following result, which is motivated by Theorem 9.11 of \cite{HL11}.
Let $\Si$ be a finite subset of $\widehat{\Gamma}_0$ satisfying the assumption (\ref{eqn::assumption}),  
and let $\bd$ be a dimension vector supported on $\Si$.
Define
\begin{equation}\label{defWd}
W^{\bd} := \bigoplus_{(i,n)\in \Si} \bC^{d_{(i, n)}}.
\end{equation}
\begin{theorem}\label{theo::iso of varieties}
There is a $G_{W^\bd}$-equivariant isomorphism of varieties
\begin{displaymath}
	\Psi : \fM_0^{\bullet}(W^\bd) \longrightarrow \rep_{\bd} \left(e_\Si \Lambda e_\Si\right).
\end{displaymath}
\end{theorem}
\demo{  We shall write for short $e$ instead of $e_\Si$, and $W$ instead of $W^\bd$.
We proceed in two steps.  We will first show that there is a closed immersion $\Psi$, in much the same way as in Proposition 9.4 of \cite{HL11}.  Then we will prove that $\Psi$ is surjective.

Notice first that $e\Lambda e = e \Lambda_\Si e$ is a finite-dimensional algebra with a unit.  
Let $\{ \alpha_k \}_{k=1}^r$ be a finite set of generators of $e\Lambda e$.  
We assume that each $\alpha_k$ lies in some $e_x \Lambda e_y$, for vertices $x$ and $y$ of $\widehat{\Gamma} \setminus \widehat{Q}$.  Then $\rep_{\bd}(e \Lambda e)$ is a closed subset of the affine space
\begin{displaymath}
	\bigoplus_{\alpha_k \in e_x\Lambda e_y} \Hom{\bC}(W_x, W_y).
\end{displaymath}
Its coordinate ring can therefore be expressed as a quotient of a polynomial ring $\bC [\theta_{i,j}^k]_{i,j,k}$, where $\theta_{i,j}^k$ should be seen as the function sending a linear map in $\Hom{\bC}(W_x, W_y)$ to the $(i,j)$th coordinate of this map written in matrix form for some fixed bases of $W_x$ and $W_y$.

Now fix some $\widehat{Q}_0$-graded vector space $V$.  Each $\alpha_k$ is an element of some $e_x \Lambda e_y$; it is thus a linear combination of paths in $\widehat{\Gamma}$ from $x$ to $y$.  For a representation in $\Lambda^{\bullet}(V,W)$, this linear combination of paths corresponds to a linear combination of composition of linear maps, each composition going from $W_x$ to $W_y$.  Let $\psi_{i,j}^k$ be the function in $\bC [\Lambda^{\bullet}(V,W)]^{G_V}$ corresponding to the $(i,j)$th coordinate of this map from $W_x$ to $W_y$ written in matrix form.  Define an algebra morphism by
\begin{eqnarray*}
 \bC [\rep_{\bd} e \Lambda e] & \longrightarrow & \bC [\Lambda^{\bullet}(V,W)]^{G_V} \\
      \theta_{i,j}^k & \longmapsto & \psi_{i,j}^k.
\end{eqnarray*}

According to \cite{Nak11} (see also \cite{LBP}, \cite{Lusztig98}), the functions $\psi_{i,j}^k$ generate $\bC [\Lambda^{\bullet}(V,W)]^{G_V}$.  We have thus defined a surjective morphism between coordinate rings, which yields a closed immersion $	\Psi_V : \fM_0^{\bullet}(V,W) \longrightarrow \rep_{\bd} e \Lambda e$.  Since $\fM_0^{\bullet}(W)$ is the union of all $\fM_0^{\bullet}(V,W)$, we get a closed immersion $	\Psi : \fM_0^{\bullet}(W) \longrightarrow \rep_{\bd} e \Lambda e$.

We now prove that $\Psi$ is surjective.  The action of $\Psi$ can be understood using the following diagram:
\begin{displaymath}
	\xymatrix{ \Lambda^{\bullet}(V,W)\ar[dr]\ar[rrr]^{\Phi_V} &  & & \rep_{\bd}(e \Lambda e) \\
	                            & \fM^{\bullet}_0(V,W) \ar[r] & \fM^{\bullet}_0(W) \ar[ru]_{\Psi}. & 
	}
\end{displaymath}
Here, $\Phi_V$ sends a representation $M$ of $\Lambda$ in $\Lambda^{\bullet}(V,W)$ to a representation of $e \Lambda e$ isomorphic to $\Hom{\Lambda_\Si}(\Lambda_\Si e, M)$.  To prove that $\Psi$ is surjective, it is thus sufficient to prove that for any representation $N$ in $\rep_{\bd}(e \Lambda e)$, there exists a $V$ and a representation  $M$ in $\Lambda^{\bullet}(V,W)$ such that $N \cong \Phi_V(M) \cong \Hom{\Lambda_\Si}(\Lambda_\Si e, M)$.  But such a $V$ and $M$ always exist thanks to Corollary \ref{coro::functor}.  Thus $\Psi$ is surjective.
Finally, the definition of $\Psi$ makes it $G_W$-equivariant.  The theorem is proved.
}

\begin{remark}
The same result was proved in \cite[Proposition 9.8]{HL11} for a restricted choice of $\Sigma$ determined by the representation theory of the quiver $Q$.
\end{remark}

\begin{example}\label{exam::A2}
Consider the quiver $Q = 1 \to 2$ of type $\bA_2$, with height function given by $\xi_1 = 1$ and $\xi_2 = 0$.  Put $W_{(1,2)} = W_{(2,-1)} = \bC$, and put all other $W$'s to zero.  We get the picture
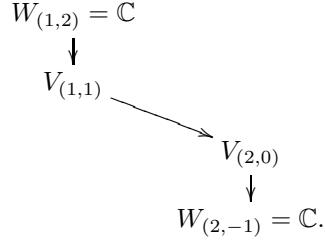
\begin{displaymath}
	\xymatrix@-1.2pc{ W_{(1,2)}=\bC\ar[d] & \\
	                  V_{(1,1)}\ar[dr] & \\
	                  & V_{(2,0)}\ar[d] \\
	                  & W_{(2,-1)}=\bC. \\
	}
\end{displaymath}
In that case, the theorem implies that $\fM_0^{\bullet}(W)$ is isomorphic to the variety of representations of $Q$ of dimension vector $(1,1)$.  
\end{example}

\begin{example}
Let $Q$ be as in Example \ref{exam::A2}.  This time, put 
\[
W_{(1,4)} = W_{(1,2)} = W_{(1,0)} = W_{(2,3)} = W_{(2,1)} = W_{(2,-1)} = \bC.
\]  
We get the picture
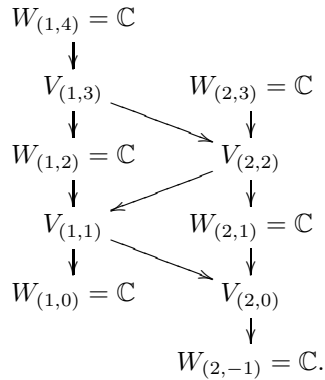
\begin{displaymath}
	\xymatrix@-1.2pc{ W_{(1,4)}=\bC\ar[d] & \\
	                  V_{(1,3)}\ar[dr]\ar[d] & W_{(2,3)}=\bC\ar[d] \\
	                  W_{(1,2)}=\bC\ar[d] & V_{(2,2)}\ar[d]\ar[dl]  \\
	                  V_{(1,1)}\ar[dr]\ar[d] & W_{(2,1)}=\bC\ar[d] \\
	                  W_{(1,0)}=\bC& V_{(2,0)}\ar[d] \\
	                  & W_{(2,-1)}=\bC. \\
	}
\end{displaymath}
The theorem states that $\fM_0^{\bullet}(W)$ is isomorphic to the variety of representations of dimension vector $(1,1,1,1,1,1)$ of the quiver
\begin{displaymath}
	\xymatrix{ 1\ar[dd]_a\ar[dddrr]^e & & \\
	                  & & 2\ar[dd]^c\ar[dddll]^g \\
	                  3\ar[dd]_b\ar[dddrr]^f & & \\
	                  & & 4\ar[dd]^d \\
	                  5 & & \\
	                  & & 6
	}
\end{displaymath}
subject to the relation $fa = -de$.
\end{example}

These examples illustrate the various forms that the algebra $e_\Sigma \Lambda e_\Sigma$ can take.  In the next section, we will see that, for a suitable choice of $W$, this algebra is the well-studied \emph{repetitive algebra}.
%-------------------------------------------------------------------------------------------------------------------------------------
\section{Orbits and strata} \label{sect3}

%................................................................................................................................
\subsection{Repetitive algebras}\label{sect3.1}
Our main source for this section is Chapter II.2 of \cite{Happel}.  Let $A$ be a finite-dimensional algebra over $\bC$.  The \emph{repetitive algebra} $\widehat{A}$ of $A$ is the infinite-dimensional algebra (without unit) defined as a vector space by
\begin{displaymath}
	\widehat{A} = \bigoplus_{n\in \bZ} A \oplus \bigoplus_{n\in \bZ} DA,
\end{displaymath}
where $D$ is the duality functor $\Hom{\bC}(?, \bC)$.  An element of $\widehat{A}$ will be denoted by $(a_n, \varphi_n)_n$, with only finitely many of the $a_n$ and $\varphi_n$ being non-zero.  The multiplication rule is given by
\begin{displaymath}
	(a_n, \varphi_n)_n \cdot (b_n, \psi_n)_n = (a_nb_n, a_{n+1}\psi_n + \varphi_n b_n)_n.
\end{displaymath}
The repetitive algebra $\widehat{A}$ can be interpreted as an ``infinite matrix algebra'' as follows:
\begin{equation*}
\widehat{A} = \left(
\begin{array}{ccccc}
\ddots &  &  & & \\
\ddots & A & & & \\
 & DA & A &  & \\
 & & DA & A  & \\
 & & & \ddots & \ddots
\end{array} \right).
\end{equation*}
Then the product in $\widehat{A}$ of two infinite matrices of this form is the usual product of matrices, except that the entries of the second lower diagonal in the product are ignored and set to zero.

As shown in \cite{Happel}, $\widehat{A}$ is a selfinjective algebra.  Thus its stable module category 
$\underline{\MOD} \widehat{A}$ is triangulated, and the suspension functor is given by the inverse of the syzygy functor $\Omega$.

\begin{example}
If $A = \bC Q$, where $Q$ is the linearly oriented quiver of type $\bA_n$
\[
 Q = 1 \to 2 \to \cdots \to n,
\]
then it is easy to see that $\widehat{A}$ is isomorphic to the path algebra of a linearly
oriented quiver $Q^{\mathrm{repet}}$ with vertex set $\bZ$ and 
the following relations:
\begin{itemize}
\item every path of length $> n$ is equal to $0$. 
\end{itemize}
A representation of $\widehat{A}$ of dimension vector $\bd=(d_i)$ is the same
as a degree 1 endomorphism $x$ of a graded vector space
\[
 W^{\bd} = \bigoplus_{i\in \bZ} \bC^{d_i}
\]
satisfying $x^{n+1} = 0$. The orbits of $\rep_\bd(\widehat{A})$ are therefore
identical to the graded nilpotent orbits  of $W$ considered by Ginzburg
and Vasserot in \cite{GV}.

\end{example}
Following J.~Schr\"oer \cite{S99}, we can describe the repetitive algebra $\widehat{\bC Q}$
of any  quiver $Q$  by a quiver $Q^{\mathrm{repet}}$ with relations as follows.
\begin{itemize}
	\item the vertices of $Q^{\mathrm{repet}}$ are labelled by $i[n]$, where $i\in Q_0$ and $n\in \bZ$;
	\item for any arrow $i \to j$ of $Q$ and any integer $n$, there is an arrow $i[n] \to j[n]$;
	\item for any maximal path $w:i \to j$ of $Q$ and any integer $n$, there is an arrow $w^*[n]:j[n] \to i[n+1]$ (these are called \emph{connecting arrows}).
\end{itemize}
Note that in \cite{S99}, the arrows $w^*[n]$ would go from $j[n]$ to $i[n-1]$; this is because our definition for the repetitive algebra, which follows \cite{Happel}, uses different conventions than those in \cite{S99}.
The relations are obtained in the following way:
\begin{itemize}
	\item a \emph{full path} is a path of the form $u[n+1]w^*[n]v[n]$, where $w=vu$ is a maximal path in $Q$.  Then any path which is not a subpath of a full path is a relation;
	\item if $w_1 = x_1vu_1$ and $w_2 = x_2vu_2$ are two maximal paths in $Q$, then the element $u_1[n+1]w_1^*[n]x_1[n] - u_2[n+1]w_2^*[n]x_2[n]$ is a relation.
\end{itemize}

\begin{example}\label{exampleA4}
Let $Q$ be the quiver of type $\bA_4$ given by
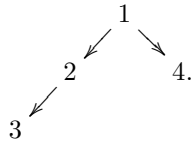
\begin{displaymath}
	\xymatrix@-1.2pc{ & & 1\ar[dr]\ar[dl] & \\
	           & 2\ar[dl] & & 4. \\
	           3 & & &
	}
\end{displaymath} 
The quiver $Q^{\mathrm{repet}}$ is:
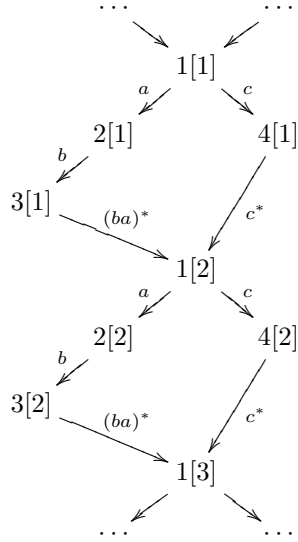
\begin{displaymath}
	\xymatrix@-1.2pc{ & \cdots\ar[dr] & & \cdots\ar[dl] \\
	                  & & 1[1]\ar[dl]_a\ar[dr]^c & \\
	                  & 2[1]\ar[dl]_b & & 4[1]\ar[ddl]^{c^*} \\
	                  3[1]\ar[drr]^{(ba)^*} & & & \\
	                  & & 1[2]\ar[dl]_a\ar[dr]^c & \\
	                  & 2[2]\ar[dl]_b & & 4[2]\ar[ddl]^{c^*} \\
	                  3[2]\ar[drr]^{(ba)^*} & & & \\
	                  & & 1[3]\ar[dl]\ar[dr] & \\
	                  & \cdots & & \cdots 
	}
\end{displaymath}
The relations are:
\[
c(ba)^* = ac^* = ba(ba)^*b =0,\qquad (ba)^*ba = c^*c.
\]
\end{example}

\begin{example}\label{exampleD4}
Let $Q$ be the quiver of type $\bD_4$ given by
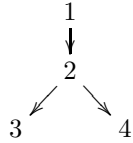
\begin{displaymath}
	\xymatrix@-1.2pc{ &  1\ar[d] & \\
	           & 2\ar[dl]\ar[dr] & &  \\
	           3 & &4 &
	}
\end{displaymath}
The quiver $Q^{\mathrm{repet}}$ is:
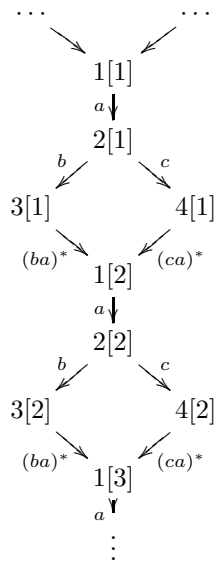
\begin{displaymath}
	\xymatrix@-1.2pc{  \cdots\ar[dr] & & \cdots\ar[dl] \\
	                  &  1[1]\ar[d]_a & \\
	                  & 2[1]\ar[dl]_b\ar[dr]^c & &  \\
	                  3[1]\ar[dr]_{(ba)^*} & &4[1]\ar[dl]^{(ca)^*} & \\
	                  &  1[2]\ar[d]_a & \\
	                  & 2[2]\ar[dl]_b\ar[dr]^c & &  \\
	                  3[2]\ar[dr]_{(ba)^*} & &4[2]\ar[dl]^{(ca)^*} & \\
	                  &  1[3]\ar[d]_a & \\
	                  & \vdots &	}
\end{displaymath} 
The relations are:
\[
ca(ba)^* = ba(ca)^*= a(ca)^*ca = a(ba)^*ba =0,\qquad (ba)^*b = (ca)^*c.
\]
\end{example}

We now recall Happel's theorem:
\begin{theorem}[Theorem II.4.9 of \cite{Happel}]\label{theo::Happel}
Assume that $A$ has finite global dimension.  Then there is an equivalence of 
triangulated categories $\cD^b(\MOD A) \rightarrow \underline{\MOD} \widehat{A}$.
\end{theorem}

In view of this theorem, the Auslander-Reiten quiver of $\MOD \widehat{A}$ is obtained by adding vertices (corresponding to the indecomposable projective-injective $\widehat{A}$-modules) to the Auslander-Reiten quiver of $\cD^b(\MOD A)$.

Assume now that $A$ is the path algebra of a quiver $Q$ of Dynkin type $\bA$, $\bD$ or $\bE$.  Then the Auslander-Reiten quiver of $\cD^b(\MOD A)$ is known to be isomorphic to $\widehat{Q}^{\mathrm{op}}$.  In order to picture the Auslander-Reiten quiver of $\MOD \widehat{A}$, we must know  which irreducible morphisms start or end in each indecomposable projective-injective module.

Let $P$ be an indecomposable projective-injective module.  Its socle and top are simple; thus any submodule and any quotient of $P$ is indecomposable.  Moreover, using for instance Proposition IV.3.5 of \cite{ASS}, we know that the only irreducible morphisms involving $P$ are the inclusion $\rad P \hookrightarrow P$ and the surjection $P \twoheadrightarrow P/\soc P$.  Therefore $P$ appears in exactly one mesh of the Auslander-Reiten quiver of $\MOD \widehat{A}$, and this mesh has the form
\begin{displaymath}
	\xymatrix@-1.2pc{  & M \ar[ddr] & \\
	            & \vdots & \\
	           \rad P \ar[uur]\ar[r]\ar[ddr] & P\ar[r] & P/ \soc P \\
	           & \vdots & \\
	           & N \ar[uur] & \\
	}
\end{displaymath}
where all modules other than $P$ are non-projective.  Thus the Auslander-Reiten
quiver of $(\MOD \widehat{A})^{\mathrm{op}}$ is a full subquiver of $\widehat{\Gamma}$. 

\begin{example}\label{exampleA4_2}
We continue Example~\ref{exampleA4}.
The Auslander-Reiten quiver of $\widehat{\bC Q}$ 
(with arrows going from bottom to top, a convention which will be useful later)
is :
\begin{displaymath}
	\xymatrix@-1.0pc{ \cdots & \cdots & \cdots & \cdots \\
	           {\begin{smallmatrix} 3[0] \end{smallmatrix}}\ar[ur] & & {\begin{smallmatrix} & 2[0] & & \\ 3[0] & & & 4[0] \\ & & 1[1] & \end{smallmatrix}}\ar[ur]\ar[ul]\ar[u] & \\
	           & {\begin{smallmatrix} 3[0] & & & 4[0] \\ & & 1[1] & \end{smallmatrix}}\ar[ur]\ar[ul] & & {\begin{smallmatrix} & 2[0] & \\ 3[0] & & \\ & & 1[1]  \end{smallmatrix}}\ar[ul] \\
	           {\begin{smallmatrix} & 4[0] \\ 1[1] & \end{smallmatrix}}\ar[ur] & & {\begin{smallmatrix} 3[0] & & \\ & & 1[1]  \end{smallmatrix}}\ar[ur]\ar[ul] & \boxed{\begin{smallmatrix} & 2[0] & \\ 3[0] & & \\ & & 1[1] \\ & 2[1] & \end{smallmatrix}}\ar[u] \\
	           \boxed{\begin{smallmatrix} & 4[0] \\ 1[1] & \\ & 4[1] \end{smallmatrix}}\ar[u]  & {\begin{smallmatrix} 1[1] \end{smallmatrix}}\ar[ur]\ar[ul] & & {\begin{smallmatrix} 3[0] & & \\ & & 1[1] \\ & 2[1] & \end{smallmatrix}}\ar[ul]\ar[u] \\
	           {\begin{smallmatrix} 1[1] & \\ & 4[1] \end{smallmatrix}}\ar[ur]\ar[u] & & {\begin{smallmatrix} & 1[1] \\ 2[1] & \end{smallmatrix}}\ar[ur]\ar[ul] & \boxed{\begin{smallmatrix} 3[0] & & \\ & & 1[1] \\ & 2[1] & \\ 3[1] & & \end{smallmatrix}}\ar[u] \\
	           & {\begin{smallmatrix} & 1[1] & \\ 2[1] & & 4[1] \end{smallmatrix}}\ar[ur]\ar[ul] & & {\begin{smallmatrix} & & 1[1] \\ & 2[1] & \\ 3[1] & & \end{smallmatrix}}\ar[ul]\ar[u] \\
	           {\begin{smallmatrix} 2[1] \end{smallmatrix}}\ar[ur] & & {\begin{smallmatrix} & & 1[1] & \\ & 2[1] & & 4[1] \\ 3[1] & & &  \end{smallmatrix}}\ar[ur]\ar[ul] & \\
	           & {\begin{smallmatrix} & 2[1] \\ 3[1] & \end{smallmatrix}}\ar[ur]\ar[ul] & \boxed{ \begin{smallmatrix} & & 1[1] & \\ & 2[1] & & 4[1] \\ 3[1] & & & \\ & & 1[2] & \end{smallmatrix}}\ar[u] & {\begin{smallmatrix} 4[1] \end{smallmatrix}}\ar[ul] \\
	           {\begin{smallmatrix} 3[1] \end{smallmatrix}}\ar[ur] &  & {\begin{smallmatrix} & 2[1] & & \\ 3[1] & & & 4[1] \\ & & 1[2] &  \end{smallmatrix}}\ar[ur]\ar[ul]\ar[u] & \\
	           \cdots & \cdots\ar[ur]\ar[ul] & \cdots & \cdots\ar[ul]
	}
\end{displaymath}
where representations are written as their composition series (each $i[r]$ represents the simple representation associated to the corresponding vertex), and projective-injective representations are written inside boxes.
\end{example}

Let $\psi$ denote a fixed embedding of the Auslander-Reiten quiver of $(\MOD \widehat{A})^{\mathrm{op}}$
as a full subquiver of $\widehat{\Gamma}$.
\begin{lemma}\label{lemm::projectives}
Let $\Si:=\psi(\proj \widehat{A})$.  Then there is an algebra isomorphism 
\[
e_\Si \Lambda e_\Si \cong \widehat{A}.
\]
\end{lemma}
\demo{The category $(\MOD \widehat{A})^{\mathrm{op}}$ is equivalent to the category generated by its Auslan\-der-Reiten quiver, together with the mesh relations.  Let $\Lambda'$ be the path algebra of this Auslander-Reiten quiver quotiented by the mesh relations.  Then $\Lambda'$ is isomorphic to the endomorphism algebra of the direct sum of all finite-dimensional indecomposable modules over $\widehat{A}$.  Thus 
\[
e\Lambda' e \cong \End{\widehat{A}}\left(\bigoplus_{P\ \mathrm{projective}} P\right)^{\mathrm{op}} \cong \widehat{A}. 
\]
The quivers of $\Lambda$ and $\Lambda'$ are the same; their relations differ only by some signs.  Consider the following vertex-preserving automorphism of $\Lambda$:
\begin{eqnarray*}
 a_{(i,n+1)} & \longmapsto & a_{(i,n+1)};\\
 b_{(i,n)} & \longmapsto & b_{(i,n)};\\
 (c,n) & \longmapsto & \begin{cases}
                      (c,n) & \text{if } n \equiv 0,1 \mod 4;\\
                      -(c,n) & \text{else};
                     \end{cases} \\
 (\overline{c},n+1) & \longmapsto & \begin{cases}
                      (\overline{c},n+1) & \text{if } n \equiv 0,1 \mod 4;\\
                      -(\overline{c},n+1) & \text{else};
                     \end{cases} \\
\end{eqnarray*}
for any $(i,n) \in \widehat{Q}_0$.  Under this automorphism, the relations of $\Lambda$ become the mesh relations; thus $\Lambda$ and $\Lambda'$ are isomorphic, and this isomorphism preserves the idempotents associated to the vertices.  So we have that $e\Lambda e \cong e \Lambda' e \cong \widehat{A}$.
}

\begin{remark}
 The isomorphism $\Lambda \cong \Lambda'$ defined above induces a $G_V$-equivariant isomorphism of varieties $\rep_\bd(\Lambda) \cong \rep_\bd(\Lambda')$, for any dimension vector $\bd$.  Thus we will often identify the representations of $\Lambda$ with those of $\Lambda'$ and conversely.
\end{remark}

\begin{corollary}\label{coro::iso of varieties}
Let $\bd$ be a finite dimension vector supported on $\psi(\proj \widehat{A})$, 
and let $W^\bd$ be the corresponding graded space defined as in 
(\ref{defWd}). 
There is a $G_W$-equivariant isomorphism of varieties 
\begin{displaymath}
	\Psi:\fM_0^{\bullet}(W^\bd) \longrightarrow \rep_{\bd}\widehat{A}.
\end{displaymath}
\end{corollary}
\demo{This follows from Lemma \ref{lemm::projectives} and Theorem \ref{theo::iso of varieties}.
}

%.....................................................................
\subsection{The Grothendieck group for repetitive algebras}

Let $\cA$ be an essentially small abelian category.  Its \emph{split Grothendieck group} 
$K_0^{\spl}(\cA)$ 
is the free abelian group generated by isomorphism classes of objects in $\cA$, subject to the relations
\begin{displaymath}
	[X] - [X\oplus Z] + [Z] = 0, \quad \textrm{for all objects $X$ and $Z$ of $\cA$.}
\end{displaymath}
The \emph{Grothendieck group of $\cA$} is the abelian group $K_0(\cA)$ obtained by adding the relations
\begin{displaymath}
	[X] - [Y] + [Z] = 0, \quad \textrm{for all exact sequences $0\rightarrow X \rightarrow Y \rightarrow Z \rightarrow 0$.}
\end{displaymath}
Let $\phi : K_0^{\spl}(\cA) \longrightarrow K_0(\cA)$ be the canonical surjection.  In \cite{Auslander84}, M.~Auslander gave a basis of $\Ker \phi$ in the case when $\cA$ is the category of finite-dimensional modules over a finite-dimensional $k$-algebra $A$ that is representation finite.  We will give a similar description of $\Ker \phi$ for the module category of a repetitive algebra of Dynkin type.  Although the proof follows the lines of \cite{Auslander84}, we include it here for completeness.

Let $\widehat{A}$ be the repetitive algebra of the path algebra of a quiver of Dynkin type.  
For any vertex $i$ of the Gabriel quiver $Q^{\mathrm{repet}}$ of $\widehat{A}$, 
let $\widehat{S}_i$ be the corresponding simple $\widehat{A}$-module, and 
let $\widehat{P}_i$ be its projective cover. 

By \cite[Lemma 2.5]{HW83}, the category $\MOD \widehat{A}$ has almost split sequences.  For any non-projective indecomposable module $M$ in $\MOD \widehat{A}$, there exists an almost split sequence (unique up to isomorphism)
\begin{displaymath}
	\xymatrix{0 \ar[r] & \tau M \ar[r] & E_M \ar[r] & M \ar[r] & 0}
\end{displaymath}
For any indecomposable module $M$, let $r_M$ be the element of $K_0^{\spl}(\MOD \widehat{A})$ defined by
\begin{equation*}
r_M = \left\{
\begin{array}{rl}
[M] - [E_M] + [\tau M] & \text{if $M$ is not projective},\\
 \phantom{} [M] - [\rad M] & \text{if $M$ is projective}.
\end{array} \right.
\end{equation*}
Define a bilinear form $h$ on $K_0^{\spl}(\MOD \widehat{A})$ by 
\begin{displaymath}
	h([M],[N]) = \dim_k \Hom{\widehat{A}}(M,N), \ \text{for any modules $M$ and $N$}.
\end{displaymath}
We gather in the next lemma the properties of $h$ that we will need.

\begin{lemma}\label{lemm::h}
\begin{enumerate}
	\item[(i)] For any indecomposable module $M$, we have $h([M], r_M) = 1$.  

	\item[(ii)] If $M$ and $N$ are indecomposable and non isomorphic, then $h([M], r_N) = 0$.
	
	\item[(iii)] For any $x \in K_0^{\spl}(\MOD \widehat{A})$, we have that
	   \begin{displaymath}
	     x = \sum_{M\ \mathrm{indec.}} h(x, r_M)[M].
     \end{displaymath}
  \item[(iv)] If $x$ and $y$ are elements of $K_0^{\spl}(\MOD \widehat{A})$ such that, for any indecomposable 
module $L$, we have that $h(x,[L]) = h(y,[L])$, then $x=y$.   The same holds if $h([L], x) = h([L], y)$ for any indecomposable $L$.
\end{enumerate}
\end{lemma}
\demo{ We first prove statement (i).  Assume that $M$ is indecomposable and is not projective.  Then there is an exact sequence 
\[
0 \to \Hom{\widehat{A}}(M, \tau M) \to \Hom{\widehat{A}}(M, E_M) \to \Hom{\widehat{A}}(M, M) \to S(M) \to 0
\]
and $h([M], r_M) = \dim S(M)$. Notice that $S(M)$ is isomorphic to the only simple $\End{\widehat{A}}(M)$-module, which is one-dimensional.  Now, if $M$ is projective, then the equality $h([M], r_M)=1$ follows from the fact that $M/\rad M$ is a one-dimensional simple module.

Statement (ii) follows from the definition of almost split sequences; indeed, if $M$ and $N$ are indecomposable and non isomorphic, with $N$ non projective, then we have an exact sequence
\begin{displaymath}
	\xymatrix{ 0 \ar[r] & \Hom{\widehat{A}}(M, \tau N) \ar[r] & \Hom{\widehat{A}}(M, E_N) \ar[r] & \Hom{\widehat{A}}(M, N) \ar[r] & 0}.
\end{displaymath}
If $N$ is projective, then any morphism from $M$ to $N$ has its image contained in $\rad N$, so that $\Hom{\widehat{A}}(M, N) = \Hom{\widehat{A}}(M, \rad N)$.  The result follows.

The equation of (iii) follows directly from (i) and (ii).  
To prove (iv), notice that if $h(x,[L]) =0$ for any indecomposable $L$, then $x = 0$ by (3).  
Thus, if $h(x, [L]) = h(y, [L])$ for any indecomposable $L$, then $x-y$ has to be zero.  The second equality is obtained by applying the first to the opposite algebra $\widehat{A}^{\mathrm{op}}$ and applying the duality functor.
}

\begin{proposition}\label{prop::groth}
Let $A = \bC Q$ for a quiver $Q$ of Dynkin type.  Let $\widehat{A}$ be the repetitive algebra of $A$.
\begin{enumerate}
	\item[(i)] For any $x\in K_0^{\spl}(\MOD \widehat{A})$, we have that
	  \begin{displaymath}
	     x = \sum_{M\ \mathrm{indec.}} h([M], x)r_M.
     \end{displaymath}
	\item[(ii)] The set $\{ r_M \ | \ M \text{ is indecomposable} \}$ is a basis for $K_0^{\spl}(\MOD \widehat{A})$.
	\item[(iii)] The set $\{ r_M \ | \ M \text{ is indecomposable and non projective} \}$ is a basis for $\Ker \phi$.
\end{enumerate}
\end{proposition}
\demo{To prove (i), notice first that the sum is finite. This follows from the fact that for a given $N$, there is a finite number of indecomposable modules $M$ such that $h([M], [N])$ is non-zero.  Thus the sum in (i) is an element of $K_0^{\spl}(\MOD \widehat{A})$.  Next, for any indecomposable module~$L$, we have that
\begin{eqnarray*}
	h\left([L], \sum_{M\ \mathrm{indec.}} h([M], x)r_M\right) & = & \sum_{M\ \mathrm{indec.}} h([M], x)h([L], r_M) \\
	                                                                & = & h([L], x)h([L], r_L) \\
	                                                                & = & h([L], x).
\end{eqnarray*}
Here, the second equality follows from Lemma \ref{lemm::h}(ii).  Applying Lemma \ref{lemm::h}(iv), we have proved (i).

It follows from (i) that the set $\{ r_M \ | \ M \text{ is indecomposable} \}$ generates the group $K_0^{\spl}(\MOD \widehat{A})$.  The fact that its elements are linearly independent follows from Lemma \ref{lemm::h}(i) and (ii): if $\sum_{M \text{ indec.}}\lambda_M r_M = 0$, then applying $h([L], ?)$ yields the equality $\lambda_L = 0$, for any indecomposable $L$.  This proves (ii).

To prove (iii), notice that $\phi(r_M) = 0$ for any non-projective indecomposable $M$, and that 
$\phi(r_{\widehat{P}_i}) = [\widehat{S}_i]$, respectively.
Thus no linear combination of the $r_{\widehat{P}_i}$ lies in the kernel of $\phi$.  This finishes the proof.
}

\begin{corollary}\label{coro::basis}
 Let $N$ be a representation of $\widehat{A}$ of dimension vector $\bd$.  Then, 
in $K_0^{\spl}(\MOD \widehat{A})$ we have
 \begin{displaymath}
  - [N] + \sum_i d_i[\widehat{S}_i]    = \sum_{M \ \mathrm{non-proj.}} \dim (\proj (\Omega M,N))\, r_M,
 \end{displaymath}
 where $\proj(A,B)$ is the space of morphisms from $A$ to $B$ that factor through a projective representation, and where $\Omega$ is the syzygy.  In particular, $N$ is completely determined up to isomorphism by the values of 
$\dim (\proj (M, N))$, as $M$ ranges over all non-projective indecomposables.
\end{corollary}
\demo{ The element $-[N] + \sum_i d_i[\widehat{S}_i]$ lies in the kernel of $\phi$. Thus, by Proposition \ref{prop::groth}, it can be written in a unique way as a linear combination of some $r_M$, for $M$ indecomposable non-projective,
and the coefficient of $r_M$ is 
\begin{displaymath}
 h([M], -[N] + \sum_i d_i[\widehat{S}_i]).
\end{displaymath}
  Now let $P^M$ be a minimal projective cover of $M$, so that we have a short exact sequence
\begin{displaymath}
 \xymatrix{ 0\ar[r] & \Omega M \ar[r] & P^M \ar[r] & M\ar[r] & 0
 }
\end{displaymath}
and an exact sequence
\begin{displaymath}
 \xymatrix{ 0\ar[r] & \Hom{}( M,N) \ar[r] & \Hom{}(P^M,N) \ar[r] & \Hom{}(\Omega M,N).
 }
\end{displaymath}
From this we get that $\proj (\Omega M, N) \cong \Hom{}(P^M,N)/\Hom{}(M,N)$.  Moreover, since $\Hom{}(P^M,N)$ only depends on the dimension vector of $N$, we get that it is isomorphic to $\Hom{}(P^M, \bigoplus_i \widehat{S}_i^{d_i})$, which is in turn isomorphic to $\Hom{}(M, \bigoplus_i \widehat{S}_i^{d_i})$, since $P^M$ is a minimal projective cover of $M$.  Thus we have that
\begin{eqnarray*}
 \dim \proj (\Omega M, N) & = & \dim \Hom{}(P^M, N) - \dim \Hom{}(M,N) \\
                          & = & \dim \Hom{}(M, \bigoplus_i \widehat{S}_i^{d_i}) - \dim \Hom{}(M,N) \\
                          & = & h\left([M],\ -[N] + \sum_i d_i[\widehat{S}_i]\right).
\end{eqnarray*}
This proves the first claim.  Since the $r_M$ form a basis of the kernel of $\phi$, a module $N'$ not isomorphic to $N$ will give rise to different values for the $h([M], [N'])$.  This proves the second claim.
}

%.....................................................................
\subsection{A bijection between strata and orbits}
Recall the algebra 
\[
\Lambda' = \End{\widehat{A}}\left(\bigoplus_{M\ \mathrm{ indec.}} M\right)^{\mathrm{op}},
\] 
where the sum is taken over all isomorphism classes $M$ of indecomposable representations of $\widehat{A}$.

\begin{lemma}\label{lemm::stable dominant}
Let $\bd$ be a finite dimension vector supported on $\psi(\proj \widehat{A})$.
Let $N$ be a representation of $\widehat{A}$  of dimension $\bd$.  
Consider the left $\Lambda'$-module
 \begin{displaymath}
   \underline{N} = \proj \left(\bigoplus_{M\ \mathrm{ indec.}} M,\ N\right).
 \end{displaymath}
Then $\underline{N}$, regarded as a point of $\Lambda^\bullet(V,W^\bd)$ for some $V$, 
is stable and dominant.
\end{lemma}
\demo{To prove that $\underline{N}$ is stable, notice that  if $L$ is an indecomposable non-projective representation of $\widehat{A}$ and if $f\in V_{\psi(L)} = \proj (L,N)$ is non-zero, then any projective cover $g:P\rightarrow L$ is such that the composition $f\circ g$ does not vanish, since $g$ is surjective.  Thus $g$, viewed as an element of $\Lambda'$, sends $f$ to a non-zero element of $W_{\psi(P)} = \proj (P,N)$.  Therefore the subrepresentation of $\underline{N}$ generated by $f$ is not supported on $V$.

To prove that $\underline{N}$ is dominant, we must prove the inequality
\begin{displaymath}
	\dim W_{(i,n)} - \dim V_{(i,n+1)} - \dim V_{(i,n-1)} + \sum_{j \sim i} \dim V_{(j,n)} \geq 0,
\end{displaymath}
for all $(i,n) \in \widehat{\Gamma}_0 \setminus \widehat{Q}_0$.  All these spaces lie in a ``mesh'' of the form
\begin{displaymath}
	\xymatrix{ V_{(i,n+1)}\ar[d]\ar[dr]\ar[drr]\ar[drrr] & & & \\
	           W_{(i,n)}\ar[d] & V_{(j_1, n)}\ar[dl] & \cdots\ar[dll] & V_{(j_r, n)}\ar[dlll] \\
	           V_{(i, n-1)} & & &
	}
\end{displaymath}
whose preimage by $\psi$ corresponds to an almost split sequence
\begin{displaymath}
 0\rightarrow \tau M \rightarrow \bigoplus_k E_k \rightarrow M \rightarrow 0.
\end{displaymath}
Notice that the middle term has a projective summand $P$ if, and only if, $(i,n) = \psi(P)$.  Thus, since $V_{\psi(L)} = \proj (L,N)$, the above inequality becomes
\begin{displaymath}
	- \dim \proj (M,N) - \dim \proj (\tau M,N) + \sum_{k} \dim \proj(E_k,N) \geq 0,
\end{displaymath}
where the counterpart of $\dim W_{(i,n)}$ is included in the sum and corresponds to the possible projective summand $P$.  Since projectives are injectives, if a morphism $\tau M \rightarrow N$ factors through a projective, then it factors through $E = \bigoplus_k E_k$. Thus we get a sequence
\begin{displaymath}
 \proj(M,N) \rightarrow \proj(E,N) \rightarrow \proj(\tau M,N)
\end{displaymath}
which is not exact, but whose first map is injective and second map in surjective, and such that their composition vanishes.  Thus 
\[
\dim \proj(E,N) \geq \dim \proj(M,N) + \proj(\tau M,N),
\] 
and $\underline{N}$ is dominant.
}
\begin{remark}
One can show that the
$G_V$-orbit of $\underline{N}$ is always closed.  This result, in a more
general context, is proved in a work in progress of Bernhard Keller and
Sarah Scherotzke. We are grateful to them for kindly explaining their result to us.
\end{remark}

\begin{theorem}\label{main_theorem}
Let $\bd$ be a finite dimension vector supported on $\psi(\proj \widehat{A})$.  Then the isomorphism
\begin{displaymath}
  \Psi: \fM_0^\bullet(W^\bd) \longrightarrow \rep_{\bd}(\widehat{A})
\end{displaymath}
of Corollary \ref{coro::iso of varieties} induces a bijection between the following sets:
\begin{enumerate}
 \item[(i)] the set of isomorphism classes of representations of $\widehat{A}$ with dimension vector~$\bd$;
 \item[(ii)] the set of dominant pairs $(V,W^\bd)$;
 \item[(iii)] the set of non-empty strata $\fM_0^{\bullet \mathrm{reg}}(V,W^\bd)$.
\end{enumerate}
\end{theorem}
\demo{ Lemma \ref{lemm::stable dominant} gives a map from set (i) to set (ii), sending a representation $N$ to the dominant pair $(V,W^\bd)$ associated to $\underline{N}$.  This map is injective since, by Corollary \ref{coro::basis}, $N$ is completely determined by the values of $\dim(\proj (L,N))$ for all non-projective indecomposable $\widehat{A}$-modules $L$, and these give the graded dimension of $V$.

Assume now that $(V,W^\bd)$ is a dominant pair.  Consider the element
\begin{displaymath}
  \sum_i d_i[\widehat{S}_i] - \sum_M \dim V_{\psi(\Omega M)}\, r_{M} \ \in K_0^{\spl}(\widehat{A}).
\end{displaymath}
If $L$ is an indecomposable non-projective representation, then the coefficient of $[L]$ in this element is given by
\begin{displaymath}
\sum_i d_i \delta_{\widehat{S}_i, L} - \dim V_{\psi(\Omega L)} - \dim V_{\psi(\Omega \tau^{-1} L)} + \sum_{\Omega L \rightarrow N} \dim V_{\psi(N)}.
\end{displaymath}
This corresponds to the almost split sequence ending in $\Omega \tau^{-1} L$. (Note that if $L$ is simple, then $\Omega L$ is the radical of a projective module, and the corresponding projective appears in the middle of this almost split sequence).  Since $(V,W^\bd)$ is dominant, this quantity is non-negative. 
On the other hand, if $P$ is an indecomposable projective, then the coefficient of $[P]$ is non-negative, since it appears only with non-negative coefficients in all the elements $r_M$.

Therefore $\sum_i d_i[\widehat{S}_i] - \sum_M (\dim V_{\psi{\Omega M}}) r_{M} = [N] $ for some representation $N$ of $\widehat{A}$ with dimension vector $\bd$.  This defines an injective map from the set (ii) to the set (i).  Since they are finite sets, they are in bijection.  Note that the two maps that we defined are inverse to each other.

Now, we know from Lemma \ref{lemm::stable dominant} that any representation $N$ gives rise to a representation $\underline{N}$ of $\Lambda'$ that is both stable and dominant.  Therefore the corresponding strata 
$\fM_0^{\bullet \mathrm{reg}}(V,W^\bd)$ is non-empty.  This gives an injective map from set (i) to set (iii).  Finally, if $\fM_0^{\bullet \mathrm{reg}}(V,W^\bd)$ is non-empty, then $(V,W^\bd)$ is dominant; this gives an injective map from set (iii) to set (ii), and we have a bijection.
}

The proof of Theorem \ref{main_theorem} gives us a description of the space $V$ in the dominant pair $(V,W)$ associated to a representation $N$ of $\widehat{A}$.

\begin{corollary}\label{coro::explicit bijection}
Let $N$ be a representation of $\widehat{A}$ of dimension vector $\bd$, and let $\fM_0^{\bullet \mathrm{reg}}(V,W^\bd)$ be the corresponding non-empty stratum.  Then the following integers are equal for any indecomposable non-projective $\widehat{A}$-module $M$:
\begin{enumerate}
	\item[(i)] $\dim V_{\psi M}$;
	\item[(ii)]	$\dim \proj (M,N)$;
	\item[(iii)] $\dim \Hom{\widehat{A}}(\Omega^{-1} M, \bigoplus_i \widehat{S}_i^{d_i}) - \dim \Hom{\widehat{A}}(\Omega^{-1}M, N)$;
	\item[(iv)] the coefficient of $r_{\Omega^{-1}M}$ in the equality
	  \begin{displaymath}
	    -[N] + \sum_i d_i [\widehat{S}_i] = \sum_{M\ \mathrm{non-proj.}} \lambda_M r_M.
    \end{displaymath}	
\end{enumerate}
\end{corollary}
\demo{ The equality of (i) and (ii) follows from the bijection given in Theorem \ref{main_theorem} and from the definition of $\underline{N}$. The equality of (ii) and (iii) was given at the end of the proof of Corollary \ref{coro::basis}.  The equality of (ii) and (iv) follows from Corollary \ref{coro::basis}.
}

\begin{remark}
We are grateful to Christof Geiss for sharing with us his observation of the equality of (i) and (iii) in Corollary \ref{coro::explicit bijection}. 
Note that by Nakajima's theory, the stratum $\fM_0^{\bullet \mathrm{reg}}(V,W^\bd)$ contains the stratum 
$\fM_0^{\bullet \mathrm{reg}}(V',W^\bd)$ in its closure if and only if $\dim V_{(i,n)} \ge \dim V'_{(i,n)}$
for every $(i,n)\in \widehat{Q}_0$.
Thus, using the equality of (i) and (iii), we find that a $\bd$-dimensional $\widehat{A}$-module $N'$
lies in the closure of the orbit of~$N$ if and only if 
$\dim \Hom{\widehat{A}}(M, N) \le \dim \Hom{\widehat{A}}(M, N')$ for every indecomposable
$\widehat{A}$-module $M$. Hence we recover that for the repetitive algebra of a Dynkin quiver,
which is locally of finite representation type, the degeneration order and the Hom order 
coincide \cite{Zwara99}.
\end{remark}
\begin{remark}
In case $W$ is supported on the image by $\psi$ of projectives associated to vertices $1[r], 2[r], \ldots, n[r]$ for some integer $r$, Theorem~\ref{main_theorem} recovers Theorem 9.11 of \cite{HL11}.  In that paper, the place of the $W$'s correspond to the place of the simple modules in $\MOD \bC Q$.  Our place for the $W$'s is the place of the projective indecomposable modules in $\MOD \widehat{Q}$.  To see that the two choices give the same pattern for the placing of the $W's$, notice that the syzygy $\Omega$ (which is an autoequivalence of the triangulated category $\underline{\MOD}\widehat{A}$) sends any simple $S$ to the radical $\rad P$ of a projective module $P$, and that $\rad P$ is the start of an almost split sequence having $P$ in its middle term.
\end{remark}

\begin{proposition}\label{proposition_multiplicities}
Let $N$ be a representation of $\widehat{A}$ of dimension vector $\bd$, and write 
\[
N = \bigoplus_{M\ \mathrm{indec.}}M^{a_M}.
\]  
Then, if $M$ is not projective, we have
\begin{displaymath}
a_{\Omega^{-1}\tau M} =	-\dim\proj(M,N) - \dim\proj(\tau M,N) + \sum_k \dim\proj(E_k, N),
\end{displaymath}
where $0\to \tau M\to \bigoplus_k E_k \to M \to 0$ is an almost split sequence.
\end{proposition}
\demo{ Three cases arise:
\begin{itemize}
	\item the almost split sequences ending in $M$ and in $\Omega^{-1} M$ have no projective summands in their middle terms;
	\item the almost split sequence ending in $M$ has a projective summand $\widehat{P}$ (this is equivalent to the fact that $M$ is isomorphic to $\widehat{P}/\soc(\widehat{P})$);
	\item the almost split sequence ending in $\Omega^{-1} M$ has a projective summand $\widehat{Q}$ (this is equivalent to the fact that $\Omega^{-1} M$ is isomorphic to $\widehat{Q}/\soc(\widehat{Q})$).
\end{itemize}
Let $\delta$ (resp. $\epsilon$) be equal to $1$ if $M$ (resp. $\Omega^{-1}M$) is isomorphic to $\widehat{P}/\soc(\widehat{P})$ (resp. $\widehat{Q}/\soc(\widehat{Q})$), and $0$ otherwise.
We can therefore write the almost split sequences as
\begin{displaymath}
	0 \to \tau M \to \widehat{P}^{\delta} \oplus \bigoplus_{E_k \textrm{ non-proj}} E_k \to M \to 0
\end{displaymath}
\begin{displaymath}
	0 \to  \Omega^{-1} \tau M \to \widehat{Q}^{\epsilon} \oplus \bigoplus_{E_k \textrm{ non-proj.}} \Omega^{-1} E_k \to \Omega^{-1} M \to 0.
\end{displaymath}
Then, using Corollary~\ref{coro::explicit bijection}, we get 
\begin{eqnarray*}
 & & \!\!\!\!\!\!\!\!\!\!\!\!\!\!\! -\dim\proj(M,N) - \dim\proj(\tau M,N) + \sum_k \dim\proj(E_k, N) \\
 &=& -\dim\proj(M\oplus \tau M,N) + \sum_{E_k \textrm{ non-proj.}} \dim\proj(E_k, N) + \delta\dim\proj(\widehat{P}, N) \\
 &=& -\dim\proj(M\oplus \tau M,N) + \sum_{E_k \textrm{ non-proj.}} \dim\proj(E_k, N) + \delta d_{\psi(\widehat{P})} \\
 &=& h\Big(-r_{\Omega^{-1}M} - \epsilon [\widehat{Q}], \quad \sum_{i}d_i \widehat{S_i} - [N]\Big) + \delta d_{\psi(\widehat{P})} \\
 &=& a_{\Omega^{-1}\tau M} - \sum_i d_i\delta_{\Omega^{-1}\tau M, \widehat{S}_i} + \delta d_{\psi(\widehat{P})} \\
 &=& a_{\Omega^{-1}\tau M} - \delta d_{\psi(\widehat{P})} + \delta d_{\psi(\widehat{P})} \\
 &=& a_{\Omega^{-1}\tau M}. 
\end{eqnarray*}
Here we have used that, since $\widehat{Q}$ is projective-injective, 
\[
h\Big([\widehat{Q}], \ \sum_{i}d_i \widehat{S_i} - [N]\Big) = 0.
\]
}

\begin{example}
We continue Example~\ref{exampleA4_2}.  Put $W_{(3,8)} = W_{(1,4)} = W_{(3,0)} = \bC$; these spaces correspond to the images of some projectives under some inclusion $\psi$ of $\MOD \widehat{A}$ into $\widehat{\Gamma}$, so that Theorem~\ref{main_theorem} applies.  
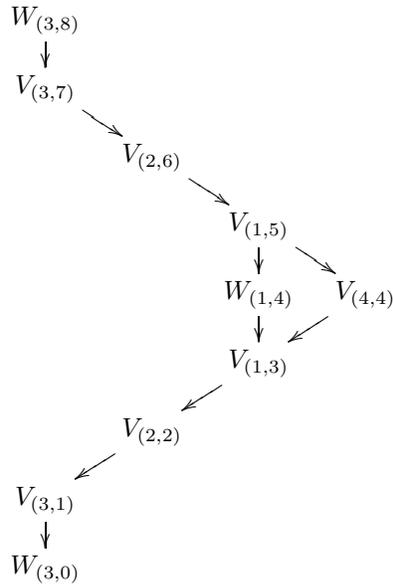
\begin{displaymath}
	\xymatrix@-1.2pc{ W_{(3,8)}\ar[d] & & & \\
	                  V_{(3,7)}\ar[dr] & & & \\
	                  & V_{(2,6)}\ar[dr] & & \\
	                  & & V_{(1,5)}\ar[d]\ar[dr] & \\
	                  & & W_{(1,4)}\ar[d] & V_{(4,4)}\ar[dl] \\
	                  & & V_{(1,3)}\ar[dl] & \\
	                  & V_{(2,2)}\ar[dl] & & \\
	                  V_{(3,1)}\ar[d] & & & \\
	                  W_{(3,0)} & & &
	}
\end{displaymath}
There are $4$ isomorphism classes of representations of $\widehat{A}$ of the corresponding dimension vector, which is supported on the following subquiver of $Q^{\mathrm{repet}}$:
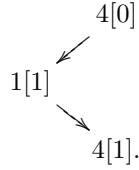
\begin{displaymath}
	\xymatrix@-1.2pc{ & 4[0]\ar[dl] \\
	                   1[1]\ar[dr] & \\
	                  & 4[1].
	}
\end{displaymath}
Corollary \ref{coro::explicit bijection} gives the dominant pair $(V,W)$ associated to each isomorphism class.  The table below describes the bijection explicitly.
\begin{table}[h]%
\begin{tabular}{|l|c|c|c|c|}
\hline
 & ${\begin{smallmatrix} 4[0] \end{smallmatrix}} \oplus {\begin{smallmatrix} 1[1] \end{smallmatrix}} \oplus {\begin{smallmatrix} 4[1] \end{smallmatrix}}$ & ${\begin{smallmatrix} & 4[0] \\ 1[1] & \end{smallmatrix}}\oplus{\begin{smallmatrix} 4[1] \end{smallmatrix}}$ & ${\begin{smallmatrix} 4[0] \end{smallmatrix}}\oplus {\begin{smallmatrix} 1[1] & \\ & 4[1] \end{smallmatrix}}$ & ${\begin{smallmatrix} & 4[0] \\ 1[1] & \\ & 4[1] \end{smallmatrix}}$ \\
\hline
$\dim V_{(3,7)}$ &0 &1 &0 &1 \\
\hline
$\dim V_{(2,6)}$ &0 &1 &0 &1 \\
\hline
$\dim V_{(1,5)}$ &0 &1 &0 &1 \\
\hline
$\dim V_{(4,4)}$ &0 &0 &0 &1 \\
\hline
$\dim V_{(1,3)}$ &0 &0 &1 &1 \\
\hline
$\dim V_{(2,2)}$ &0 &0 &1 &1 \\
\hline
$\dim V_{(3,1)}$ &0 &0 &1 &1 \\
\hline
\end{tabular}

\label{}
\end{table}
\end{example}

\begin{example}
The bijection of Theorem \ref{main_theorem} does not hold for a general $W$.  Indeed, let $Q$ be the quiver $1\to 2$ of type $\bA_2$.  Put $W_{(1,3)} = W_{(1,1)} = W_{(2,0)} = \bC$.  
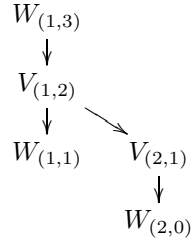
\begin{displaymath}
	\xymatrix@-1.2pc{ W_{(1,3)}\ar[d] & \\
	                  V_{(1,2)}\ar[d]\ar[dr] & \\
	                  W_{(1,1)} & V_{(2,1)}\ar[d] \\
	                  & W_{(2,0)}
	}
\end{displaymath}
Then Theorem \ref{theo::iso of varieties} gives an isomorphism between $\fM_0^{\bullet \mathrm{reg}}(W)$ and the variety of representations of dimension vector $(1,1,1)$ of the quiver
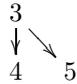
\begin{displaymath}
\xymatrix@-1.2pc{ 3\ar[d]\ar[dr] & \\ 4 & 5}	
\end{displaymath}
Then there are three dominant pairs for $W^{\bd}$, given by 
\begin{displaymath}
	(\dim V_{(1,2)}, \dim V_{(2,1)}) \in \{(0,0), (1,0), (1,1)\},
\end{displaymath}
while there are four orbits for this dimension vector, namely those of
\begin{displaymath}
	{\begin{smallmatrix} 3 \end{smallmatrix}} \oplus {\begin{smallmatrix} 4 \end{smallmatrix}} \oplus {\begin{smallmatrix} 5 \end{smallmatrix}}, \quad {\begin{smallmatrix} 3 \\4 \end{smallmatrix}} \oplus {\begin{smallmatrix} 5 \end{smallmatrix}}, \quad {\begin{smallmatrix} 3 & \\ & 5 \end{smallmatrix}} \oplus {\begin{smallmatrix} 4 \end{smallmatrix}}, \quad {\begin{smallmatrix} 3 & \\ 4 & 5 \end{smallmatrix}}.
\end{displaymath}
\end{example}

%-----------------------------------------------------------------------
\section{Composition multiplicities for quantum loop algebras} \label{sect4}

\subsection{Quantum loop algebras and $q$-characters}

Let $\g$ be a complex simple Lie algebra of type $\bA$, $\bD$ or $\bE$.
Let $Q$ be an arbitrary orientation of the Dynkin diagram of $\g$,
that is, a Dynkin quiver of the same type as $\g$.
Let $L\g$ be the loop algebra attached to $\g$, and let $U_q(L\g)$ be the
associated quantum enveloping algebra. 
We assume that the deformation parameter $q\in\bC^*$ is not a root of unity.

By \cite{FR}, every finite-dimensional $U_q(L\g)$-module $M$
(of type 1) has a \emph{$q$-character} $\chi_q(M)$, which is a 
Laurent polynomial in
\[
\cY:=\bZ\left[Y_{i,a}^{\pm1} \mid i \in Q_0,\, a\in \bC^*\right] 
\]
encoding the dimensions of its loop weight spaces. 
These $q$-characters generate a commutative subalgebra of $\cY$
isomorphic to the Grothen\-dieck ring of the category of  
finite-dimensional irreducible $U_q(L\g)$-modules.
In other words, finite-dimensional simple $U_q(L\g)$-modules are 
characterized by their $q$-characters up to isomorphism.

As in many other representation theories, the $q$-characters of the simple modules 
do not have a simple description and they are calculated via the $q$-characters
of more accessible modules called \emph{standard modules}.
As in \cite{HL11}, we will label simple modules $L(m)$ and standard modules
$M(m)$ by the unique monomial $m$ of their $q$-characters corresponding to their 
highest weight vectors. These highest weight monomials $m$ are called 
\emph{dominant} because they only involve nonnegative powers 
of the variables $Y_{i,a}$.
 
Hence, a basic question is to calculate for every pair $(m,m')$ of dominant monomials
the composition multiplicities
$\zeta_{m,m'}$ appearing in the expansions
\[
\chi_q(M(m')) = \sum_{m} \zeta_{m,m'}\, \chi_q(L(m)). 
\]

\subsection{Nakajima varieties and composition multiplicities}

By a classical reduction procedure (see \eg \cite{HL10}), it is enough to determine the multiplicities
$\zeta_{m,m'}$ corresponding to pairs $(m,m')$ of dominant monomials in the variables of the set
\[
\cZ := \left\{ Y_{i,{q^n}} \mid (i,n) \in \widehat{\Gamma}_0 \setminus \widehat{Q}_0\right\}, 
\]
where $\widehat{\Gamma}_0$ and $\widehat{Q}_0$ have been defined in \S\ref{sect2.1}.
The modules $L(m)$ where $m$ ranges over all dominant monomials in the variables
of $\cZ$ generate a tensor subcategory $\cC_\bZ$ of the category of finite-dimensional
$U_q(L\g)$-modules.

From now on, we shall only consider monomials in the variables of $\cZ$,
and write for short $Y_{i,n}$ instead of $Y_{i,q^n}$.
In particular, we define for $(i,n)\in \widehat{Q}_0$ the Laurent monomial
\[
A_{i,n} := Y_{i,n-1} Y_{i,n+1} \prod_{j\sim i} Y_{j,n}^{-1}. 
\]
Given a finite dimension vector $\bd$ on $\widehat{\Gamma}$ and the two corresponding 
graded vector spaces $V$ and $W$ as in (\ref{eqVW}), we define two monomials
\[
A^V := \prod_{(i,n)\in \widehat{Q}_0} A_{i,n}^{-d_{i,n}},
\qquad
Y^W := \prod_{(i,n)\in \widehat{\Gamma}_0 \setminus \widehat{Q}_0} Y_{i,n}^{d_{i,n}}. 
\]
It is easy to check that the pair $(V,W)$ is dominant in the sense of \S\ref{sect2.2}
if and only if the monomial $m_\bd := Y^WA^V$ is dominant, that is, every variable 
$Y_{i,n}$ has a non-negative exponent in $m_\bd$.

Given a dominant pair $(V,W)$ such that $\fM^{\bullet\, {\rm reg}}_0(V,W) \not = \emptyset$,
denote by $IC_W(V)$ the intersection cohomology complex of the closure of 
the stratum $\fM^{\bullet\, {\rm reg}}_0(V,W)$. Let $\cH^i(IC_W(V))$ be its $i$th
cohomology sheaf, and $\cH^i(IC_W(V))_{V'}$ be the stalk of this sheaf at a point
of $\fM^{\bullet\, {\rm reg}}_0(V',W)$.
We can now state one of the main results of Nakajima's geometric approach to the
representation theory of $U_q(L\g)$.

\begin{theorem}{\rm\cite[\S8]{N}}\label{NakajimaTh}
Let $m$ and $m'$ be two dominant monomials in the variables of $\cZ$. 
The multiplicity of the simple module $L(m)$ as a composition factor of the 
standard module $M(m')$ is given by
\begin{equation}\label{Nakajima_form}
\zeta_{m,m'} = 
\sum_{i\ge 0}  \dim\cH^i(IC_W(V))_{V'}, 
\end{equation}
for any pair of strata $\fM^{\bullet\, {\rm reg}}_0(V,W)$ and $\fM^{\bullet\, {\rm reg}}_0(V',W)$ 
such that $m = Y^WA^V$ and $m'=Y^WA^{V'}$.
\end{theorem}

\subsection{Parametrization of irreducible modules}

Let $N$ be a representation of $\widehat{A}$ of dimension $\bd$. 
By Theorem~\ref{main_theorem}, $N$ corresponds to a unique dominant
pair $(V,W^{\bd})$. This allows us to attach to $N$ a dominant monomial
\[
 m_N := Y^{W^\bd} A^V
\]
in the variables of $\cZ$.

Let $N'$ be another $\widehat{A}$-module of dimension $\bd'$.
Suppose that $N \simeq N' \oplus P$, where $P$ is a projective $\widehat{A}$-module.
Let $(V',W^{\bd'})$ and $(V'',W'')$ be the dominant pairs associated with $N'$
and $P$ by Theorem~\ref{main_theorem}.
Then we have the following isomorphisms of graded vector spaces:
\[
W^{\bd} \simeq W^{\bd'} \oplus W'',\qquad V \simeq V' \oplus V''.  
\]
The second isomorphism follows from the fact that, by Corollary~\ref{coro::explicit bijection}, 
the dimensions of the graded components of $V$ are equal to the dimensions of the spaces $\proj(M,N)$ for 
the indecomposable non-projective $\widehat{A}$-modules $M$, and these dimensions are additive
with respect to~$N$. 
Moreover, since $P$ is projective-injective, 
\[
\dim\proj(M,P) = \dim\Hom{\widehat{A}}(M,P),
\] 
which is additive on the almost split sequences 
$\tau M \to E \to M$ of the proof of Lemma~\ref{lemm::stable dominant}. 
It follows that 
\begin{displaymath}
	\dim W''_{(i,n)} - \dim V''_{(i,n+1)} - \dim V''_{(i,n-1)} + \sum_{j \sim i} \dim V''_{(j,n)} = 0,
\end{displaymath}
for every $(i,n)\in \widehat{\Gamma}_0\setminus \widehat{Q_0}$.
Therefore $Y^{W''}A^{V''} = 1$, and 
\[
m_N= Y^{W^\bd} A^V = Y^{W^{\bd'}} A^{V'} = m_{N'}.
\]
Hence the monomial $m_N$ 
depends only on the isomorphism
class of $N$ in the stable category $\underline{\MOD} \widehat{A}$. 

So we can assume that $N$ has no projective indecomposable summand
and write
\[
N = \bigoplus_{M\ \mathrm{indec. non proj.}}M^{a_M}.
\]  
Let us write the dominant monomial $m_N$ as
\[
m_N = \prod_{(i,n)\in \widehat{\Gamma}_0 \setminus \widehat{Q}_0} Y_{i,n}^{b_{i,n}}, 
\]
and let us express the exponents $b_{i,n}$ in terms of the multiplicities $a_M$.
By definition of $m_N$, we have
\begin{displaymath}
b_{i,n} =	\dim W^{\bd}_{(i,n)} - \dim V_{(i,n+1)} - \dim V_{(i,n-1)} + \sum_{j \sim i} \dim V_{(j,n)}.
\end{displaymath}
Now, arguing as in the proof of Lemma~\ref{lemm::stable dominant},
we get 
\[
b_{i,n} = - \dim \proj (M,N) - \dim \proj (\tau M,N) + \sum_{k} \dim \proj(E_k,N)
\]
where the almost split sequence
\begin{displaymath}
 0\rightarrow \tau M \rightarrow \bigoplus_k E_k \rightarrow M \rightarrow 0.
\end{displaymath}
corresponds under $\psi$ to the mesh containing the spaces $W^\bd_{i,n}$, $V_{i,n+1}$,
$V_{i,n-1}$, and $V_{j,n}$ for $j\sim i$. In other words, $M = \psi^{-1}(i,n+1)$.
Now, applying Proposition~\ref{proposition_multiplicities},
we get 
\[
b_{i,n} = a_{\Omega^{-1}\tau \psi^{-1}(i,n+1)}. 
\]
In particular we see that the indecomposable non-projective 
$\widehat{A}$-modules $N$ have their corresponding monomials
$m_N$ equal to the variables $Y_{i,n}$ of $\cZ$.
%Hence we see that if we denote by $\theta$ the bijection from 
%$\widehat{\Gamma}_0 \setminus \widehat{Q}_0$ to the set of isoclasses 
%of indecomposable non-projective $\widehat{A}$-modules defined by
%\[
%\theta(i,n) = \Omega^{-1}\tau \psi^{-1}(i,n-1), 
%\]
Thus, summarizing the above discussion, we have proved:
\begin{proposition}\label{parametrization}
The map $N \mapsto L(m_N)$ induces a bijection between the isomorphism 
classes of objects in the stable category $\underline{\MOD} \widehat{A}$
and the simple objects of $\cC_\bZ$. 
\end{proposition}

\subsection{Composition multiplicities}

Let $M(m')$ be a standard module in the subcategory $\cC_\bZ$, and let
$N'$ be a module over the repetitive algebra $\widehat{A}$
such that $m' = m_{N'}$. By Proposition~\ref{parametrization},
$N'$ exits and its isoclass is completely determined by $m'$,
up to projective summands. Let $\bd$ be the dimension vector
of $N'$.
We can now state the main result of this section.

\begin{theorem}
The multiplicities of the composition factors of 
$M(m')$ are equal to 
\[
\zeta_{m_N,m'} = 
\sum_{i\ge 0}  \dim\cH^i(IC(\cO_N))_{N'}. 
\]
Here $\cO_N$ denotes the orbit of a module $N$ in $\rep_{\bd}\widehat{A}$
containing $N'$ in its closure, 
$IC(\cO_N)$ is the intersection cohomology complex of the closure of $\cO_N$,
and $\cH^i(IC(\cO_N))_{N'}$ is the stalk of the $i$th cohomology sheaf of
$IC(\cO_N)$ at $N'$.
\end{theorem}

\demo{
%By Lemma~\ref{lem_monomial}, we can write $m' = Y^{W^\bd}A^{V'}$
%for some dimension vector $\bd$ supported on $\psi(\mathrm{proj}\widehat{A})$ 
%and some $\widehat{Q}_0$-graded space $V'$.
Let $(V',W^{\bd})$ be the dominant pair corresponding to $N'$. 
It is well known that if $\zeta_{m,m'}\not = 0$
we must have $m = m'A^{V''}$ for 
some $\widehat{Q}_0$-graded space $V''$. 
Hence, putting $V:=V'\oplus V''$,
we have $m= Y^{W^\bd}A^{V}=m_N$ for some $\widehat{A}$-module $N$
of dimension $\bd$.
By Corollary~\ref{coro::iso of varieties}, the Nakajima variety
$\fM_0^{\bullet}(W^\bd)$ is isomorphic to $\rep_\bd(\widehat A)$.
%the variety of representations
%of dimension $\bd$ of the repetitive algebra $\widehat{A}$ of $Q$.
By Theorem~\ref{main_theorem}, the strata
$\fM_0^{\bullet \mathrm{reg}}(V,W^\bd)$ 
and $\fM_0^{\bullet \mathrm{reg}}(V',W^\bd)$
are non-empty and isomorphic to the orbits $\cO_N$ and $\cO_{N'}$.
The theorem then follows from Theorem~\ref{NakajimaTh}.
}

%--------------------------------------------------------------------------
\def\cprime{$'$}
\providecommand{\bysame}{\leavevmode\hbox to3em{\hrulefill}\thinspace}
\providecommand{\MR}{\relax\ifhmode\unskip\space\fi MR }
% \MRhref is called by the amsart/book/proc definition of \MR.
\providecommand{\MRhref}[2]{%
  \href{http://www.ams.org/mathscinet-getitem?mr=#1}{#2}
}
\providecommand{\href}[2]{#2}

%\bibliographystyle{amsplain} 
%\bibliography{cluster_new}

\end{document}